\magnification=\magstep1
\baselineskip=1.3\baselineskip

\def\wh{\widehat}

\def\qed{{\hfill $\square$ \bigskip}}
\def\sqr#1#2{{\vcenter{\vbox{\hrule height.#2pt
        \hbox{\vrule width.#2pt height#1pt \kern#1pt
           \vrule width.#2pt}
        \hrule height.#2pt}}}}
\def\square{\mathchoice\sqr56\sqr56\sqr{2.1}3\sqr{1.5}3}

\centerline{\bf COALESCENCE OF SKEW BROWNIAN MOTIONS}
\vskip0.7truein
\centerline{{\bf Martin Barlow},$^1$ 
\footnote{}{1. Research partially supported by an NSERC (Canada) grant.} 
{\bf Krzysztof Burdzy},$^2$
\footnote{}{2. Research partially supported by NSF grant DMS-9700721.}
{\bf Haya Kaspi}$^3$ and {\bf Avi Mandelbaum}$^3$
\footnote{}{3. Research
partially supported by the Fund for the Promotion of Research at the Technion.}}
\vskip1truein

\noindent
The purpose of this short note is to prove
almost sure coalescence of two skew Brownian motions starting
from different initial points, assuming that they are driven
by the same Brownian motion. The result is very simple
but we would like to record it in print as it has already become
the foundation of a research project of Burdzy and Chen (1999).
Our theorem is a by-product of an investigation
of variably skewed Brownian motion, see Barlow et al. (1999).

Suppose that $B_t$ is the standard Brownian motion
with $B_0=0$
and consider the equation
$$X^x_t = x + B_t + \beta L_t^x, \qquad t\ge 0,  \eqno(1) $$
where $X^x_t$ satisfies the initial condition $X^x_0=x$.
Here $\beta$ is a fixed number in $[-1,1]$
and $L^x_t $ is the symmetric local time of $X^x_t$ at $0$.
Harrison and Shepp (1981) proved that (1) has a unique strong 
solution, which is skew Brownian motion. One way to define
skew Brownian motion in the case $\beta\geq 0$
is to start with a standard Brownian motion $B'_t$ 
and flip every excursion of $B'_t$ below 0 to the positive
side with probability $\beta$, independent of what
happens to other excursions. See It\^o and McKean (1965) or
Walsh (1978) for more information.

\bigskip\noindent
{\bf Theorem}. {\sl
If $X^x_t$ and $X^y_t$ are solutions
of (1) with the same $\beta\in[-1,1]\setminus\{0\}$,
relative to the same Brownian motion $B_t$,
then $X^x_t = X^y_t$ for some $t<\infty$, a.s.
}

\bigskip\noindent
{\bf Proof}.
For simplicity assume that $\beta >0$ and $0=x< y$.
Let $\wh L^0_t = \beta L^0_t$, $\wh L^y_t = y+\beta L^y_t$, and 
$$\eqalign{
T_0 &= 0, \cr
S_k &= \inf \{t> T_k: -B_t = \wh L^y_t\}, \ \ k\geq 0, \cr
T_k &= \inf \{t> S_{k-1}: -B_t = \wh L^0_t\}, \ \ k\geq 1, \cr
W_k &= {\wh L^y_{S_{k-1}} - \wh L^0 _{S_{k-1}} \over
\wh L^y_{T_{k-1}} - \wh L^0 _{T_{k-1}} }, \ \ k\geq 1, \cr
V_k &= {\wh L^y_{T_{k}} - \wh L^0 _{T_{k}} \over
\wh L^y_{S_{k-1}} - \wh L^0 _{S_{k-1}} }, \ \ k\geq 1, \cr
M_k &= \wh L^y_{T_{k}} - \wh L^0 _{T_{k}}, \ \ k\geq 0. }$$

We will first find the distributions of $W_k$'s and $V_k$'s
using excursion theory. Recall the fundamentals of excursion
theory for the standard Brownian motion from, e.g.,
Karatzas and Shreve (1991). The Brownian excursions from $0$
form a Poisson point process whose clock can be identified
with the local time of Brownian motion at $0$.
The intensity of excursions on the positive side of $0$
whose height is greater than $h$ is equal to $1/(2h)$.

The stopping time $S_0$ may be described as the first time
when an excursion of $-B_t-\wh L^0_t$ above $0$ hits the level
$y-\wh L^0_t$. 
These excursions can be identified with the excursions
of the skew Brownian motion $X^0_t$ below $0$.
They form a Poisson point
process $\cal P$ similar to the Poisson point process of excursions
of the standard Brownian motion from 0. The intensity
of $\cal P$-excursions above $0$
with height greater than $h$ is equal to
$(1 - \beta)/(2h)$. Note the extra factor $1-\beta$ as compared
to the analogous formula for the excursions of the
standard Brownian motion. The factor can be explained using
the excursion flipping construction of skew Brownian
motion mentioned in the introduction---in a sense,
the fraction of excursions flipped to the other side
is equal to $\beta/2$.
When the clock $L^0_t$ for the Poisson point process $\cal P$
takes a value $u$ then the instanteneous intensity
of excursions with height greater than $y-\wh L^0_t$
is equal to $(1 - \beta)/(2(y- \beta u))$.
We have $\wh L^y_{S_{0}} - \wh L^0 _{S_{0}} < a$ if 
no $\cal P$-excursion with height greater than $y - \wh L^0_t$
occurs before the time $s$ when $y - \wh L^0_s = a$,
i.e., when $L^0_s = (y-a)/\beta$.
Thus excursion theory enables us to write the probability
of this event using Poisson
probabilities as follows,
$$P(\wh L^y_{S_{0}} - \wh L^0 _{S_{0}} < a)
=
\exp\left( - \int_0^{(y-a)/\beta}{1 - \beta\over 2(y- \beta u)} du\right)
= \left({a \over y}\right)^{(1-\beta)/(2\beta)}.$$
Recall that 
$\wh L^y_{T_{0}} - \wh L^0 _{T_{0}} = y$.
We have
$$
P(W_1 y < a)=
P(W_1 (\wh L^y_{T_{0}} - \wh L^0 _{T_{0}}) < a)
= P(\wh L^y_{S_{0}} - \wh L^0 _{S_{0}} < a)
= (a/y)^{(1-\beta)/(2\beta)}.$$
By changing the variable we obtain for $w \in (0,1)$,
$$P(W_1 < w) = w ^{(1-\beta)/(2\beta)}.$$
By the strong Markov property, $P(W_k < w) = w ^{(1-\beta)/(2\beta)}$
for $w \in (0,1)$ and every $k\geq 1$.

A totally analogous argument shows that
$P(V_k > v) = v ^{-(1+\beta)/(2\beta)}$
for $v\geq 1$ and $k\geq 1$.

Note that, by the strong Markov property, all
random variables $V_k,W_k, k\geq 1$, are jointly independent.

Next we will show that the process
$M_k$ is a martingale and converges to 0.
First, note that $M_k = M_{k-1} W_kV_k$.
It is elementary to check that $EW_k = (1-\beta)/(1+\beta)$
and $EV_k = (1+\beta)/(1-\beta)$.
By the joint independence of $W_k$'s and $V_k$'s,
$$E(M_k \mid M_{k-1}, M_{k-2}, \dots)
= M_{k-1} EW_kEV_k = M_{k-1},$$
which shows that $M_k$ is a martingale.
As a positive martingale, the process $M_k$ must converge
with probability 1 to a random variable $M_\infty$.
Since for every $k$, $M_k$ is the product of $M_{k-1}$ and
an independent random variable $W_kV_k$, the limit
$M_\infty$ can take only the values $0$ or $\infty$.
By Fatou's Lemma, $EM_\infty \leq EM_0 = y$, so $M_\infty =0$
a.s.

On every interval $[T_k, S_k]$ the process
$\wh L^y_t - \wh L^0_t$ is non-increasing
but it is non-decreasing on intervals of the form
$[S_k, T_{k+1}]$. Thus
$$\sup_{t\in[T_k, T_{k+1}]} \wh L^y_t - \wh L^0_t
\leq \max (M_k, M_{k+1}).$$
In view of convergence of $M_k$ to 0, we must
have a.s.~convergence of $\wh L^y_t - \wh L^0_t$
to $0$ when $t\to\infty$.
It remains to show that the convergence does not take
an infinite amount of time.

Let $T_\infty = \lim_{k\rightarrow \infty} T_k$. In
view of the remarks in the last paragraph, it is not hard to see that
the value of
$\wh L^0_{T_\infty}$ is bounded by
$\sum_{k=1}^\infty M_k$. Since $\wh L^0 _\infty = \infty$, it will suffice
to show that
$\sum_{k=1}^\infty M_k < \infty$
in order to conclude that $T_\infty < \infty$. 
We have for $k\geq 1$,
$$M_k = y \prod_{j=1}^{k} W_j V_{j}.$$ We can write 
$$
y \prod_{j=1}^{k} W_j V_{j}
= \exp \left (
\log y + \sum_{j=1}^{k}
\left[ \log W_j + \log V_{j} \right ] \right) .$$ 
One can directly check
that the distribution of $-\log W_j$ is exponential with mean
$2\beta/(1-\beta)$, while the distribution of $\log V_j$ is exponential with
mean
$2\beta/(1+\beta)$. Thus, $E( \log W_j + \log V_{j+1} ) < 0$.
It follows that for some $a >0$,
we eventually have
$$\sum_{j=1}^{k}
\left[ \log W_j + \log V_{j} \right ]
\leq -a k.$$
Hence, for some random $c_1$ and all $k$
we have $M_k \leq c_1 e^{-ak}$ and so $\sum_{k=1}^\infty M_k <
\infty$, a.s. \qed

\bigskip
\centerline{\bf References}

\item{[1]} Barlow, M., Burdzy, K., Kaspi, H. and Mandelbaum, A.
(1999), Variably skewed Brownian motion (preprint)
\item{[2]} Burdzy, K. and Chen, Z.-Q. (1999) Local time flow
related to skew Brownian motion (preprint)
\item{[3]} Harrison, J.M.\ and Shepp, L.A. (1981), On skew Brownian
motion, {\it
Ann.\ Probab}. {\bf 9} (2), 309--313.
\item{[4]} It\^o, K.\ and McKean, H.P. (1965), {\it Diffusion Processes
and Their Sample Paths}, Springer, New York.
\item{[5]} Karatzas, I.\ and Shreve, S.E. (1991), {\it Brownian Motion
and Stochastic
Calculus}, 2nd Edition, Springer Verlag, New York. 
\item{[6]} Walsh, J.B.
(1978), A diffusion with discontinuous local time, {\it Temps Locaux
Asterisque}, {\bf 52--53}, 37--45. 

\vskip1truein

\noindent
Martin Barlow: University of British Columbia,
Vancouver, BC V6T 1Z2, Canada \hfil\break {\it barlow@math.ubc.ca}

\noindent
Krzysztof Burdzy: University of Washington, Seattle, WA 98195-4350,
USA \hfil\break {\it burdzy@math.washington.edu}

\noindent
Haya Kaspi and Avi Mandelbaum: Technion Institute, Haifa, 32000, Israel
\hfil\break
{\it iehaya@tx.technion.ac.il, avim@tx.technion.ac.il}

\bye